\documentstyle[12pt]{article}

  \newcommand{\const}{\rm const}

  \newcommand{\vraisup}{\rm vraisup}

\begin{document}

   \begin{center}

{\bf Extremal case of  parabolic  differential equations having  discontinuous unbounded coefficients.} \par

\vspace{4mm}

{\it Existence of fundamental solution for an initial Cauchy problem. Parametrix method}\par

\vspace{4mm}

 {\bf M.R.Formica, \ E.Ostrovsky, \ L.Sirota }

\end{center}

\vspace{6mm}

 \ Universit\`{a} degli Studi di Napoli Parthenope, via Generale Parisi 13, Palazzo Pacanowsky, 80132,
Napoli, Italy. \\

e-mail: mara.formica@uniparthenope.it \\

\vspace{3mm}

 \ Israel, Bar-Ilan University, department of Mathematic and Statistics, 59200, \\

\vspace{3mm}

   e-mail:eugostrovsky@list.ru\\

   e-mail:sirota3@bezeqint.net\\

\vspace{3mm}

\begin{center}

 {\bf Abstract} \\

\end{center}

\vspace{4mm}

 \ We prove in this short report the existence of a fundamental solution  (F.S.) for the Cauchy initial boundary problem
 on the whole space for the parabolic differential equation having at origin the point of non-integrable unbounded
  discontinuity for coefficient before a first order derivative. \par
  \ We give also the non-asymptotic rapidly decreasing at infinity estimate for these function. \par
  \ We extend the classical parametrix method offered by E.E.Levi. \par

\vspace{4mm}

\begin{center}

 {\bf Key words and phrases} \\

\end{center}

\vspace{4mm}

 \ Parabolic differential equation (PDE), parametrix method, fundamental solution (F.S.), degenerate diffusion, random variable (r.v.), random process (r.p.),
 Gamma, Beta-functions; Bessel's and variable Bessel's random process,
 modified Bessel's function, reflection at the origin, singular heat equation, subgaussian random variables and processes, Brownian motion,
 Markov random process, transfer density of probability for diffusion random processes,  Mittag-Leffler function, non-asymptotic estimate,
 Chapman-Kolmogorov equation, Neuman series, H\"older's continuity, Volterra's integral equation, initial Cauchy problem.\\

\vspace{5mm}

\section{Definitions. Notations. Statement of problem. }

\vspace{5mm}

 \ We consider in this report the Partial Differential Equation (PDE) of (uniform) parabolic type (singular heat equation)

\begin{equation} \label{equation1}
\frac{\partial u}{\partial t} = \frac{1}{2} \frac{\partial^2 u}{\partial x^2} + \frac{ 1 + 2b(t,x)}{x} \ \frac{\partial u}{\partial x},
\end{equation}
$ \ u = u(t,x), \ $
where $ \ t \ge 0, \ x \in R_+, \   $ and  the coefficient $ \ b = b(t,x) \ $ is continuous bounded function
such that in general case $ \ b(t,0) \ne 0, \ $  common with the condition of (homogeneous)  reflection at the origin

\begin{equation} \label{refl origin}
\lim_{x \to 0+} x^{1 - 2 b(t,x)} \ u'_x(t,x) = 0
\end{equation}

and initial value problem (Cauchy statement)

\begin{equation} \label{initcond}
\lim_{t \to s+} u(t,x) = f(x), \ u = u(t,x) = u[f](t,x), \ s \ge 0,
\end{equation}
 where $ \ f(\cdot) \  $ is certain measurable function satisfying some growth condition at  $ \ x \to \infty; \ $  for definiteness one can suppose its continuity and boundedness:
$ \ \vraisup_x |f(x)| < \infty. \ $ \par
 \ Briefly, $ \ L_{t,x}u = L[b]_{t,x}u= 0, \ $ where $ \ L[b]_{t,x}  \ $ is a linear parabolic partial differential operator of a parabolic type
 with singular not integrable coefficients   of the form

\begin{equation} \label{operator L}
L[b]_{t,x} = L_{t,x} \stackrel{def}{=} \frac{\partial }{\partial t} - \frac{1}{2} \frac{\partial^2 }{\partial x^2} - \frac{1 + 2b(t,x)}{x} \ \frac{\partial} {\partial x}
\end{equation}
in the domain $ \ t,x > 0. \ $ \par

 \vspace{4mm}

 \ Let us recall the classical definition. \par

\vspace{4mm}

{\bf Definition 1.1.} The (measurable) function $ \ p = p(t,x,s,y), \ 0 < s < t; \ (x,y) \in (R_+)^2 \ $ as ordinary is said to be a fundamental solution (F.S.) for
the equation (\ref{equation1}) common with the initial condition (\ref{initcond}) and the reflection condition  (\ref{refl origin}), iff
 for all the fixed values $ \ (s,y), \ 0 < s < t, \ x,y \in R_+ \ $ it satisfies the equation (\ref{equation1} ) as well  as  the initial condition (\ref{initcond}) and the
 reflection condition  (\ref{refl origin}) and for all the bounded continuous function $ \ f = f(x) \ $

\begin{equation} \label{funsol}
\lim_{t \to s+} \int_{R_+} p(t,x, s,y) \ f(y) \ dy  = f(x), \ x \in R_+,
\end{equation}
if of course there exists.  \par
 \ In this case the F.S. $ \ p(t,x,s,y) \ $ may be interpreted  as the
 transfer density of probability for diffusion random non-homogeneous in the time as well as in the space Markov process, see e.g.
\cite{Bally}, \cite{Corielli}.  As a consequence:

$$
p(t,x, s,y)  \ge 0; \ \int_0^{\infty} p(t,x, s,y) \ dy = 1
$$
and it satisfies also the famous Chapman-Kolmogorov equation

\begin{equation} \label{Chapman Kolm}
p(t,x,s,y) = \int_0^{\infty} p(t,x,v,z) \ p(v,z,s,y) \ dz, \ 0 < s < v < t.
\end{equation}

 \vspace{4mm}

  \ On the other words, the (unique) solution $ \ u = u(t,x) \ $  of the Cauchy problem  (\ref{equation1}), (\ref{initcond}),  (\ref{refl origin})
 in the whole space $ \ R_+ \ $ may be written as follows

 \begin{equation} \label{solution}
 u(t,x) = \int_{R_+} p(t,x,s,y) \ f(y) \ dy.
 \end{equation}

  \vspace{3mm}

  \ For instance, the famous Gauss function

\begin{equation} \label{Zfun}
Z(t,x,s,y) = Z(t-s, x-y) := (2 \pi (t-s))^{-1/2} \ \exp \left\{ \  - \frac{1}{2} \ \frac{(x-y)^2}{(t-s)}  \ \right\}
\end{equation}
is the F.S. for the ordinary heat equation only with initial condition in the whole space $ x \in R \ $

$$
\frac{\partial u}{\partial t}  =  \frac{1}{2} \ \frac{\partial^2 u}{\partial x^2}.
$$

\vspace{4mm}

 \ The case of the reflected Brownian motion will be considered further. \par

\vspace{4mm}

 \ Another example, important for us. Consider the following infinitesimal operator for diffusion random Markov process  $ \ \xi =\xi(t) \ $

$$
L_a = \frac{1}{2} \frac{\partial ^2}{\partial x^2} + \frac{1 + 2a}{x} \frac{\partial}{\partial x}, \ a = \const \in (-1,0), \ x \in R_+,
$$
with reflection at the origin

$$
f \in D(L_a) \Rightarrow  \left\{ f \in C^2(0,\infty), \  \lim_{x \to 0+} x^{1 - 2 a} f'(x) = 0  \ \right\}
$$
and corresponding PDE

$$
\frac{\partial u}{\partial t} = L_a \ u
$$
with ordinary initial condition. \par

\vspace{3mm}

 \ If $ \ a \ge 0 \ $ then the point $ \ x = 0 \ $ is unattainable; when $ \ a \le - 1 \ $ it is absorbing point, see \cite{Molchanov}. {\it Therefore
 we consider in this report only  the case } $ \ - 1 < a < 0 \ $  {\it and correspondingly}

 $$
 -1 < b(t,x) < 0.
 $$

\vspace{3mm}

  \ The fundamental solution has a form $ \ p(t,x,s,y) =  p_a(t,x,s,y) \stackrel{def}{=} \  $

\begin{equation} \label{FSconstcoef}
  (t-s)^{-1} \ x^{-a} \ y^{a + 1} \ \exp \left\{ \ - \frac{x^2 + y^2}{2(t-s)} \ \right\} \ I_a \left( \frac{xy}{t-s}  \right),
\end{equation}
$ \ x,y > 0, \ t >s, \ $ where $ \ I_a(\cdot) \ $ denotes an usually {\it modified} Bessel's function of order $ \ a: \ $

\begin{equation} \label{Mod Bessel}
I_a(z) \stackrel{def}{=} \sum_{k=0}^{\infty} \frac{(z/2)^{a + 2k}}{k! \ \Gamma(a + k + 1)},
\end{equation}
see \cite{Molchanov},  \cite{Nikeghbali}. \par

 \ The case $ \  a = -1/2 \ $ correspondent to the ordinary reflected at the origin Brownian motion:

$$
p_{-1/2}(t,x,s,y) = (2 \pi (t-s))^{-1/2} \ \left\{ \exp \left[ \frac{(x + y)^2}{2(t-s)}  \right] + \exp \left[ \frac{(x - y)^2}{2(t-s)}  \right]   \right\},
$$
 $ \ x,y > 0, \ 0 < s < t. \ $ \par

 \vspace{4mm}

 \ Note by the way as long as

$$
 \int_0^{\infty} p(t,x,s,y) dy = 1,
$$

\vspace{4mm}

$$
\int_0^{\infty} y^{a + 1} \ e^{-y^2/2t } \ I_a(xy/t) \ dy = t  \ x^a v \ e^{x^2/2t }, \ x,y,t > 0, \ a \in (-1,0),
$$
or equally for at the same values of parameters

\begin{equation} \label{Int Bes}
\int_0^{\infty} z^{a+1} \ e^{ -w z^2/2  } \ I_a(z) \ dz = w^{-a - 1} e^{1/2w}, \ w > 0.
\end{equation}

 \ Further, it follows from the Chapman-Kolmogorov's equation (\ref{Chapman Kolm}) that

$$
\int_0^{\infty} (t-v)^{-1} \ x^{-a} \ z^{a+1}  \ \exp \left\{ - \frac{x^2 + z^2}{2(t-v)}  \right\} \ I_a \left(  \frac{xz}{t-v} \right) \times
$$

$$
(v-s)^{-1}  \ z^{-a} \ y^{a+1} \ \exp \left\{ -\frac{z^2 + y^2}{2(v-s)} \right\}  \ I_a \left( \frac{zy}{v-s}   \right) \ dz =
$$

\begin{equation} \label{Conv estim}
(t-s)^{-1} \ x^{-a} \ y^{a+1} \ \exp \left\{ - \frac{x^2 + y^2}{2(t-s)}  \right\} \ I_a \left( \frac{xy}{t-s}  \right), 0 < s < v < t < \infty,
\end{equation}
or equally

$$
\int_0^{\infty} t^{-1}  \ z  \ \exp \left\{ - \frac{x^2 + z^2}{2t}  \right\} \ I_a \left(  \frac{xz}{t} \right) \times
$$

$$
s^{-1} \ \exp \left\{ -\frac{z^2 + y^2}{2s} \right\}  \ I_a \left( \frac{zy}{s}   \right) \ dz =
$$

\begin{equation} \label{Kol Chap2}
(t+s)^{-1} \ \exp \left\{ - \frac{x^2 + y^2}{2(t+s)}  \right\} \ I_a \left( \frac{xy}{t+s}  \right), 0 < s, t < \infty.
\end{equation}

\vspace{4mm}

 \ {\it  We will emphasize that the condition }  $ \   a \in (-1,0)  \ $ {\it is very important for the last fourth  relations. } \par

\vspace{5mm}

 {\bf  Our claim in this preprint is to prove the existence of F.S. for our equation  (\ref{equation1}) with initial condition (\ref{initcond})
 and reflection condition (\ref{refl origin}) under suitable conditions.  } \par

\vspace{4mm}

 \ We give also the exact bilateral estimates for this fundamental solution. \par

\vspace{5mm}

 \ The existence of F.S. for {\it strictly}  parabolic PDE with smooth coefficients, for instance, bounded and belonging to the H\"older space, is  explained in many works, e.g.
 \cite{Bally}, \cite{Dressel2}, \cite{DeckKruse}, \cite{Formica M}, \cite{Friedman}, \cite{Gagliardo}, \cite{IlKalashOleynik}, \cite{Ladyzhenskaya}, \cite{Lax}, \cite{Levi1}, \cite{Levi2}.
 At the same problem for parabolic systems is considered in \cite{Eydel},  \cite{Slobodetskiy}. The case of parabolic equations with discontinuous coefficients is investigated in
\cite{DeckKruse}, \cite{Ilyin}, \cite{Kamynin},   \cite{Krasnosel}, \cite{Samarskiy}, \cite{Sobolevskiy} etc.\par
 \ The recent results about the parabolic equations with singular coefficients may be found in articles
\cite{Byun}, \cite{Paul Eric}, \cite{Voit}, \cite{Zhang} and so one. \par

\vspace{4mm}

 \ We will essentially follow the interest article  of  Thomas Deck and Susanne Kruse  \cite{DeckKruse}, in which was considered the case  of an equation

\begin{equation} \label{equation3}
\frac{\partial u}{\partial t} = \frac{1}{2} \frac{\partial^2 u}{\partial x^2} + b_1(t,x) \ \frac{\partial u}{\partial x},
\end{equation}
in which the coefficient $ \ b_1 = b_1(t,x) \ $ is H\"older continuous and may be unbounded at $ \ x \to \pm \infty \ $ in the following sense

$$
|b_1(t,x) | \le C \  \left( \ 1 + |x|^{\beta} \ \right), \ C,\beta = \const < \infty.
$$

 \ The uniqueness  of the F.S. for parabolic PDE  is proved, e.g. in  \cite{Dressel1},  \cite{Dressel2}, \cite{Gagliardo}, \cite{Kamynin}, \cite{Ladyzhenskaya}, \cite{Samarskiy}. \par

\vspace{4mm}

 \ In the probability theory, more precisely, in the theory of diffusion random processes the F.S. $ \ p(t,x,s,y) \ $ represent the density of transition of probability.
 The equations of the form
 (\ref{equation1}) appears in particular in \cite{Molchanov}, \cite{Paul Eric}; see also \cite{Bally}, \cite{Corielli}. \par

\vspace{5mm}

\section{Main result. Existence of the F.S. }

\vspace{5mm}

 \ {\bf Theorem 2.1.}  Suppose in the equation (\ref{equation1}) the coefficient $ \ b = b(t,x)  \ $ is   H\"older continuous:

$$
\exists \alpha = \const \in (0,1], \ H = \const \in (0,\infty) \ \Rightarrow
$$

\begin{equation} \label{Holder}
 | b(t,x) - b(s,y)| \le H [ |t-s| + |x-y|   ]^{\alpha}.
\end{equation}

 \ Suppose also that this function is bilateral bounded as follows

\begin{equation} \label{restriction}
 \beta \stackrel{def}{=} \inf_{t,x} b(t,x) > - 1,
\end{equation}

\begin{equation} \label{restriction1}
 \beta_+ \stackrel{def}{=} \sup_{t,x} b(t,x) < 0.
\end{equation}

 \ Then there exists a F.S. $ \ p(t,x,s,y) \ $ for the equation (\ref{equation1}), common with initial and reflection conditions,
which satisfies in addition the following estimates  (\ref{estim for p}), \ (\ref{lower estim for p}). \par

\vspace{5mm}

 \ {\bf Proof.} \par

 \ {\bf 1.} \ We will use the classical parametrix method, see \cite{Levi1}, \cite{Levi2}, \cite{Dressel2}, \cite{DeckKruse}, \cite{Friedman},
 \cite{Gagliardo}, \cite{IlKalashOleynik}, \cite{Ladyzhenskaya}, \cite{Lax},  etc., applying as the initial approximation the function $ \ p_{b(s,y)}(t,x,s,y) \ $
 instead the well - known approximation $ \ Z(t,x,s,y). \ $ \par
  \ Note that the initial approximation of the form $ \ Z(t,x,s,y). \ $ was used  also in the recent article \cite{Formica M} for the parabolic PDE
  having discontinuous but local integrable coefficient

\begin{equation} \label{equation11 gam}
\frac{\partial u}{\partial t} = \frac{1}{2} \frac{\partial^2 u}{\partial x^2} + \frac{ 1 + 2b(t,x)}{|x|^{\gamma}} \ \frac{\partial u}{\partial x},
\end{equation}
 $ \gamma = \const \in (0,1), \  x \in R. \ $\par

 \vspace{4mm}

  \  We will follow the article \cite{IlKalashOleynik}, down to notations. \par
 \  Introduce for this purpose the following functions.

$$
w = w(t,x,\tau,\xi) =  w_{\tau,\xi}(t,x,\tau,\xi) := p_{b(\tau,\xi)}(t,x,\tau,\xi),
$$

$$
w^{(2)} = w^{(2)}(t,x,\tau,\xi) = w_{t,x}(t,x,\tau,\xi) := p_{b(t,x)}(t,x,\tau,\xi),
$$
where (here  and throughout whole report)

$$
\ t,x,\tau,\xi > 0, \tau < t.
$$

 \ The F.S. $ \ p(t,x,\tau, \xi) \ $  for source equation \ (\ref{equation1}) \ may be found from the representation

\begin{equation}  \label{form for p}
p(t,x,\tau,\xi) :=  w_{\tau,\xi}(t,x,\tau,\xi) + \int_{\tau}^t d \theta \ \int_0^{\infty} \ w_{\theta,\zeta}(t,x,\theta,\zeta) \cdot \Phi(\theta,\zeta,\tau,\xi) \ d\zeta,
\end{equation}
see \cite{IlKalashOleynik}, relation 4.18, where the function $ \ \Phi = \Phi(\theta,\zeta,\tau,\xi) \ $ may be found from the following linear integral equation

\begin{equation}  \label{equat for p1}
\Phi(t,x, \tau,\xi)= L_{t,x} \left( \ w_{\tau,\xi}(t,x,\tau,\xi) \ \right) +
\end{equation}

\begin{equation}  \label{equat for p2}
\int_{\tau}^t \ d \theta \ \int_0^{\infty} L_{t,x}\left(w_{\theta,\zeta}(t,x,\theta,\zeta) \right) \cdot
\Phi(\theta,\zeta,\tau,\xi) \ d \zeta.
\end{equation}
see \cite{IlKalashOleynik}, relation 4.54. \par

\vspace{5mm}

\ {\bf 2.} Denote for convenience

$$
K = K(t,x,\tau,\xi) :=L_{t,x} \left( \ w_{\tau,\xi}(t,x,\tau,\xi) \ \right), \ t,x,\tau,\xi > 0;
$$
and define  for brevity following  the authors of \cite{DeckKruse} the two-dimensional (generalized) convolution

\begin{equation} \label{conv}
f*g(t,x,s,y) \stackrel{def}{=} \int_s^t \ dr \int_0^{\infty} f(t,x,r,z) g(r,z,s,y) \ dz,
\end{equation}
then  the equation (\ref{equat for p2}) takes the form $ \ \Phi(t,x,s,y) = K(t,x,s,y) + \ $

\begin{equation} \label{equat K1}
 \int_s^t dr \int_0^{\infty} dz K(t,x,r,z) \cdot \Phi(r,z,s,y), 0 < s \le t < \infty, \ x,y > 0.
\end{equation}
or briefly

$$
 \ \Phi = K + K*\Phi.
$$

  \ The expression for  the kernel $ \  K(t,x,s,y),  \ $ as well as for the "right-hand" side of (\ref{equat K1}) has a form

\begin{equation} \label{expr for K}
K(t,x,s,y) = \frac{b(t,x) - b(s,y)}{x} \cdot \frac{\partial p_{b(s,y)}(t,x,s,y)}{\partial x}.
\end{equation}

\vspace{5mm}

 \ {\bf 3.} \ The last two equations are well-known Volterra's equations. Its (unique) solution may be represented by the form

 \begin{equation} \label{series Phi}
\Phi(t,x,s,y) = K(t,x,s,y)  +\sum_{m=1}^{\infty} \Phi_m(t,x,s,y),
 \end{equation}
 where

 $$
 \Phi_1(t,x,s,y) = K*K(t,x,s,y)
 $$
and recursively

$$
\Phi_{m+1}(t,x,s,y)  = \Phi_m *K(t,x,s,y), \ m = 1,2,\ldots.
$$
 \ Define formally $ \  \Phi_0(t,x,s,y) := K(t,x,s,y);   \ $  so that

\begin{equation} \label{series K}
\Phi(t,x,s,y) = \sum_{m=0}^{\infty} \Phi_m(t,x,s,y),
 \end{equation}

\vspace{4mm}

 \ It remains to ground the convergence  of the series (\ref{series Phi}) or (\ref{series K}) and made some quantitative estimates.\par

\vspace{5mm}

 \ {\bf 4.} Note first of all that  it follows from the condition (\ref{Holder}) the estimate

\begin{equation} \label{first estim}
|K(t,x,s,y)| \le H \ \frac{(|x-y| + |t-s|)^{\alpha}}{x} \  \left| \ \frac{\partial p_{b(s,y)}(t,x,s,y)}{\partial x} \ \right|.
\end{equation}
 \ Further, we will use the following well known properties of the modified Bessel's functions: the power series in (\ref{Mod Bessel}) and the following
 asymptotical expression

\begin{equation} \label{As Bess}
z \to \infty \ \Rightarrow I_a(z) \sim (2 \pi z)^{-1/2} \ \exp(z).
\end{equation}

\ We will use also some estimations relative (modified) Bessel's and gamma functions from the works \cite{Abramowitz}, \cite{Abramowitz Bes}, \cite{Arfken}, \cite{Barkat},
\cite{Bateman}, \cite{Gradshteyn}, \cite{Kazeminia}.

\vspace{4mm}

 \ Let $ \ \delta \ $ be arbitrary fixed number from an open interval $ \ (0,1): \ \delta \in (0,1). \  $
  We get synthesizing both the estimates for Bessel's  functions (\ref{As Bess}) and (\ref{Mod Bessel}) analogously to the work \cite{IlKalashOleynik}

$$
|K(t,x,s,y)| \le C(\beta,\beta_+) \ (t-s)^{-1 + \alpha/2} \times
$$

\begin{equation}\label{K estim}
 \frac{H}{\sqrt{e \ \delta}} \ p_{\beta}(t, (1 - \delta)x, s, (1 - \delta)y), 0 < s < t; \ x,y > 0.
\end{equation}

\vspace{5mm}

\ {\bf 5.} We have after simple calculations using convolution identity (\ref{Conv estim})

$$
|\Phi_1(t,x,s,y)| \le C^2(\beta,\beta_+) \ (1 -\delta)^{-1} \ (e \ \delta)^{-1} \ \frac{H^2 \ \Gamma^2(\alpha/2)}{\Gamma(\alpha)}  \times
$$

$$
(t-s)^{\alpha - 1} \ p_{\beta}(t, \ (1 - \delta)x, \ s,  \ (1 - \delta)y),
$$
where as ordinary $ \ \Gamma(\cdot) \ $  denotes the Gamma function. \par

\vspace{5mm}

\ {\bf 6.} We find by means of an induction

$$
|\Phi_m(t,x,s,y)| \le  C^m(\beta,\beta_+) \ H^{m+1} \ (1 - \delta)^{-m} \ (e \ \delta)^{-m} \ \frac{\Gamma^{m+1}(\alpha/2)}{\Gamma((m+1) \alpha/2) } \times
$$

\begin{equation} \label{Phi m estim}
(t-s)^{(m+1) \alpha/2 -1} \ p_{\beta}(t,(1 - \delta)x, s, (1 - \delta)y), \ m = 0,1,2,\ldots.
\end{equation}
 \ Therefore, for all the values $ \ \delta \in (0,1) \ \Rightarrow \  $

\begin{equation} \label{series for Phi 1}
|\Phi(t,x,s,y)| \le H \ (t-s)^{-1 + \alpha/2}  \ \Gamma(\alpha/2) \ p_{\beta}(t,x,s,y) \times
\end{equation}

\begin{equation} \label{series for Phi 2}
\sum_{m=0}^{\infty} \frac{H^m \ C^m(\beta,\beta_+) \ \Gamma^m(\alpha/2) \ (1 - \delta)^{-m} \ (e \delta)^{-m} \ (t-s)^{m \ \alpha/2} }{ \Gamma((m+1)\alpha/2)}.
\end{equation}

\vspace{4mm}

 \ Notice that the series in (\ref{series for Phi 2}) convergent, following the F.S. $ \ p(t,x,s,y)  \ $  for source equation there exists. \par

\vspace{4mm}

 \ Let us estimate the obtained before function $ \ \Phi(t,x,s,y). \ $
 Recall  for this purpose briefly the classical definition of the so - called (generalized) Mittag-Leffler functions  \ \cite{Mittag}, \cite{Mittag2}, \cite{Shukla}

\vspace{4mm}

$$
E_{A,B}(z) \stackrel{def}{=} \sum_{k=0}^{\infty} \frac{z^k}{\Gamma(A k + B)}.
$$
\ Here  $ \ A = \const > 0, \ B \ge 0. \ $   We need only  further the following particular case on these functions

\begin{equation} \label{Mitt}
g_{\alpha}(z) \stackrel{def}{=} \sum_{m=0}^{\infty} \frac{z^m}{\Gamma(m\alpha/2 + \alpha/2)} = E_{\alpha/2,\alpha/2}(z), \ \alpha \in (0,1].
\end{equation}.

 \ We need to use as a values  $ \ z \ $ in (\ref{Mitt}) only real non - negative  values $ \ z \in [0,\infty), \ $ but the expression in  (\ref{Mitt})
is also defined for all the complex values $ \ z. \ $ \par

 \ The asymptotic as well as non-asymptotic behavior of these functions as $ \ |z| \to \infty \ $ is investigated in  \cite{Shukla}, see also \cite{Evgrafov}.
 Namely, as $ \ z \to \infty, \ z \in R_+ \ $

$$
E_{A,B}(z) \sim A^{-1} \ \exp \left( \ z^{1/A} \   \right).
$$
 \ As a particular  case

\begin{equation} \label{as Mittag}
g_{\alpha}(z) \sim \frac{2}{\alpha} \ \exp \left( \ z^{2/\alpha} \ \right), \ z \to \infty.
\end{equation}

 \ Put   now

$$
v = v_{\delta}(t-s) := H \ C(\beta,\beta_1) \ \Gamma(\alpha/2) \ (1 - \delta)^{-1} (e \delta)^{-1} \ (t-s)^{\alpha/2},
$$
then

$$
\tilde{v} : = v_{1/2}(t-s) = H \ C(\beta,\beta_1) \ \Gamma(\alpha/2) \ 4 \ e^{-1}  \ (t-s)^{\alpha/2},
$$

\begin{equation} \label{Phi Mittag}
|\Phi(t,x,s,y)| \le H \ (t-s)^{-1 + \alpha/2}  \ \Gamma(\alpha/2) \ p_{\beta}(t,x,s,y) \cdot g_{\alpha}(v_{\delta}(t-s)).
\end{equation}

 \ Of course

\begin{equation} \label{Phi Mittag  inf}
|\Phi(t,x,s,y)| \le \inf_{\delta \in (0,1)} H \ (t-s)^{-1 + \alpha/2}  \ \Gamma(\alpha/2) \ p_{\beta}(t,x,s,y) \cdot g_{\alpha}(v_{\delta}(t-s))
\end{equation}
and

\begin{equation} \label{Phi Mittag half}
|\Phi(t,x,s,y)| \le H \ (t-s)^{-1 + \alpha/2}  \ \Gamma(\alpha/2) \ p_{\beta}(t,x,s,y) \cdot g_{\alpha}(v_{1/2}(t-s)).
\end{equation}

\vspace{5mm}

\ {\bf 7.} Finally, it remains to apply the representation (\ref{form for p}) for the fundamental solution $ \ p(t,x,\tau,\xi) \ $.
We estimate

$$
\forall \delta \in (0,1) \ \Rightarrow    p(t,x,\tau,\xi) \le  p^{(\delta)}(t,x,\tau,\xi),
$$
where

$$
p^{(\delta)}(t,x,\tau,\xi) \stackrel{def}{=}  w_{\tau,\xi}(t,x,\tau,\xi) + \int_{\tau}^t d \theta \ \int_0^{\infty} \ w_{\theta,\zeta}(t,x,\theta,\zeta) \cdot \Phi(\theta,\zeta,\tau,\xi) \ d\zeta=
$$

$$
p_{b(\tau,\xi)}(t,x,\tau,\xi) + \int_{\tau}^t d \theta \ \int_0^{\infty} \ p_{b(\theta,\zeta)}(t,x,\theta,\zeta) \cdot \Phi(\theta,\zeta,\tau,\xi) \ d\zeta \le
$$

\begin{equation}  \label{estim for p}
p_{b(\tau,\xi)}(t,x,\tau,\xi) + C_2(\beta, \ \beta_+)  H \ (t-s)^{-1 + \alpha/2}  \ \Gamma(\alpha/2) \ p_{\beta}(t,x,\tau,\xi) \cdot g_{\alpha}(v_{\delta}(t-s)).
\end{equation}

\vspace{3mm}

\ Of course,

\begin{equation} \label{inf delta}
 p(t,x,\tau,\xi) \le \overline{p}(t,x,\tau,\xi) \stackrel{def}{=} \inf_{\delta \in (0,1)} p^{(\delta)}(t,x,\tau,\xi),
\end{equation}
and following

\begin{equation} \label{inf delta half}
 p(t,x,\tau,\xi) \le  p^{(1/2)}(t,x,\tau,\xi).
\end{equation}

\vspace{4mm}

 \ On the other hand,  $ \ \forall \delta \in (0,1) \ \Rightarrow p(t,x,\tau,\xi)  \ge \ $

\begin{equation}  \label{lower estim for p}
p_{b(\tau,\xi)}(t,x,\tau,\xi) - C_3(\beta, \ \beta_+)  H \ (t-s)^{-1 + \alpha/2}  \ \Gamma(\alpha/2) \ p_{\beta}(t,x,\tau,\xi) \cdot g_{\alpha}(v_{\delta}(t-s))
\end{equation}

\ Recall that

 $$
  p_{b(\tau,\xi)}(t,x,\tau,\xi) = (t-\tau)^{-1} \ x^{-b(\tau,\xi)} \ \xi^{b(\tau,\xi)+1} \ \exp \left\{ \ - \frac{x^2 + \xi^2}{2(t-\tau)} \ \right\} \times
 $$

$$
   I_{b(\tau,\xi)} \left( \frac{x \xi}{t-\tau}  \right), \  x,\xi > 0, \ t > \tau > 0.
$$

\vspace{4mm}

 \ One can  establish analogously from the relation  \ (\ref{form for p}) \ also the estimations for the suitable derivatives
 $ \ \partial p(t,x,s,y)/\partial t, \ \partial p(t,x,s,y)/\partial x, \ \partial^2 p(t,x,s,y)/\partial x^2 \ $ etc. to ground that the obtained function
solved the source problem.  \par

\vspace{5mm}

 \section{Solving of several partial differential equations.}

\vspace{5mm}

 {\bf 1.}  Let us return to the source problem (\ref{equation1}), (\ref{refl origin}), and (\ref{initcond}). \par

 \vspace{4mm}

 \  {\bf Proposition 3.1.} \ We suppose in the sequel that the coefficient $ \ b = b(t,x) \ $ satisfies all  the  conditions of theorem 2.1.
Suppose also that the initial function $ \ f = f(x), \ x \ge 0 \ $ is  continuous  and such that

\begin{equation} \label{cond f}
\exists \Delta \in (0,1) \ \forall x \ge 0  \Rightarrow \   |f(x)| \le \exp \left((1 - \Delta)x^2/2)    \right).
\end{equation}

 \ Then the (unique) solution of the considered problem has a form

\begin{equation} \label{sol initial}
u(t,x) = \int _0^{\infty} p(t,x,s,y) \ f(y) \ dy.
\end{equation}

\vspace{5mm}

\ {\bf Proof.} Namely, the convergence of the integral in (\ref{sol initial}) common with its derivatives follows immediately from  (\ref{cond f}) and from the estimate
 (\ref{estim for p}), if we choose for instance $ \ \delta = \Delta/2. \ $ \par

\vspace{5mm}

{\bf 2.} \  Let us consider now the non - homogeneous equation of the form

\begin{equation} \label{equation non hom}
\frac{\partial u}{\partial t} = \frac{1}{2} \frac{\partial^2 u}{\partial x^2} + \frac{ 1 + 2b(t,x)}{x} \ \frac{\partial u}{\partial x} + g(t,x),
\end{equation}
$ \ u = u(t,x), \ $ where $ \ t \ge 0, \ x \in R_+, \   $
 under reflection restriction

\begin{equation} \label{refl origin again}
\lim_{x \to 0+} x^{1 - 2 b(t,x)} \ u'_x(t,x) = 0
\end{equation}
and zero  initial value

\begin{equation} \label{hom initcond}
\lim_{t \to s+} u(t,x) = 0.
\end{equation}

\vspace{4mm}

 \ One can state analogously the following assertion. \par

 \vspace{4mm}

 \  {\bf Proposition 3.2.} \ We suppose that the coefficient $ \ b = b(t,x) \ $ again satisfies all  the  conditions of theorem 2.1.
Suppose also that the  right - hand side  function $ \ g = g(t,x), \ x \ge 0 \ $ is  continuous  and such that

\begin{equation} \label{cond f}
\exists \Delta \in (0,1) \ \forall x \ge 0  \Rightarrow \   |g(t,x)| \le \exp \left((1 - \Delta)x^2/2)    \right).
\end{equation}

 \ Then the (unique) solution of the considered problem has a form

\begin{equation} \label{sol init}
u(t,x) = \int_s^t d \theta \int _0^{\infty} p(t,x,\theta,y) \ g(\theta,y) \ dy.
\end{equation}

\vspace{5mm}

 \section{Variable Bessel's random processes.}

\vspace{5mm}

 \ Recall,  \ see  \cite{Bogus}, \cite{McKean}, \cite{Shiga}, \ that the Markov diffusion separable random  process $ \ \xi = \xi(t), \ t, \ \xi(t) \ge 0 \ $
 having an infinitesimal operator (generator) of the form

$$
L_a \ u = \frac{1}{2} \frac{\partial ^2}{\partial x^2} + \frac{1 + 2a}{x} \ \frac{\partial}{\partial x}, \ a = \const \in (-1,0), \ x \in R_+, \ a = \const
$$
with or without reflection at the origin, is named as (ordinary) {\it Bessel's} process.\par
 \ More information  about these processes may be found in the articles \cite{Bogus}, \cite{McKean}, \cite{Shiga} etc. A very interest applications of these processes
 is represented in a recent article \cite{Voit}. \par

 \ An interest example. Let $ \ w = w(t), \ t \ge 0, \ w(t) \in R^d, \ d = 1,2,\ldots \ $ be a  multivariate Brownian motion. Define for Euclidean norm $ \ ||x||, \ x \in R^d \ $
the random process $ \ \eta = \eta(t) \ $ as the (Euclidean) norm of the Brownian motion  $  \ \eta(t) = ||w(t)||. \ $ Then $ \ \eta(t) \ $ is the Bessel's random process with
correspondent generator

\begin{equation} \label{Bessel n}
M_n \ u = \frac{1}{2} \frac{\partial ^2 \ u}{\partial x^2} + \frac{n-1}{2 \ x} \ \frac{\partial u}{\partial x}.
\end{equation}

\vspace{4mm}

\ {\bf Definition 4.1.} The random diffusion process $ \ \theta(t), \ \theta(s) = x_0 = \const > 0, \ x \in R_+, \ $ having appeared before (\ref{operator L}) the infinitesimal operator

\begin{equation} \label{inf operator Qb}
Q[b]_{t,x} \ u \stackrel{def}{=} \frac{\partial^2 \ u }{\partial x^2} + \frac{1 + 2b(t,x)}{x} \ \frac{\partial \ u} {\partial x}
\end{equation}
in the domain $ \ t,x > 0 \ $ \ will be named as {\it variable Bessel's random process.} \par

\vspace{4mm}

 \ We can and will assume that the coefficient $ \ b = b(t,x) \ $ satisfies all  the  conditions of theorem 2.1., therefore $ \ \theta(t) \ $ there exists and obeys the transfer density
of certain random Markov's diffusion process $ \ \xi = \xi(t), \ t, \xi(t) > 0; \ $ defined on the suitable sufficiently rich probability space $ \ (\Omega = \{\omega\}, B, {\bf P}),\ $
and investigated before. \par

\ We intend to prove in this section that this r.p.  $ \ \xi = \xi(t) \ $ is {\it subgaussian} and continuous almost everywhere (i.e., with probability one).\par

\vspace{3mm}

 \ Let us recall at first some facts from the theory of subgaussian r.v. and r.p. More information may be found in a works
 \cite{Buldygin Koz 1}, \cite{Buldygin Koz 2}, \cite{Kahane1}, \cite{Kahane2}, \cite{Kos Os}, \cite{Ledoux Tal}, \cite{Ostrovsky1} and so one. \par

 \  The numerical valued r.v. $ \ X = X(\omega) \ $ is said to be subgaussian, if
 there exists positive finite constants $ \ C, v, \ $ such that

\begin{equation} \label{def subg}
\forall u \ge 0 \ \Rightarrow {\bf P} (|X| > u) \le C e^{-v \ u^2}.
\end{equation}

 \ These r.v. forms a so-called Birnbaum-Orlicz Banach complete space  relative the classical Luxemburg  norm

$$
||X||_{\psi(2)} \stackrel{def}{=} \inf \left\{ \ s, s > 0, \ {\bf E} \exp(|X/s|^2) \le 2 \ \right\}.
$$
 \ Equivalent (up to multiplicative constant) definition:

$$
|||X|||_{\psi(2)} \stackrel{def}{=} \sup_{p \ge 1} \left[ \ \frac{ \left[{\bf E}|X|^p \right]^{1/p}}{\sqrt{p}} \ \right].
$$

\vspace{3mm}

\ By  definition the random process $ \ X = X(t), \ t \ge 0 \ $ is named subgaussian, iff for all the values $ \ t \ge 0 \ $

$$
||X(t)||_{\psi(2)} < \infty.
$$

\ Alike, the random process $ \ X= X(t), \ t \ge 0 \ $ is named uniform subgaussian, iff

$$
\sup_{t \ge 0} ||X(t)||_{\psi(2)} < \infty.
$$

\ Analogously, the random process $ \ X= X(t), \ t \ge 0 \ $ is named local subgaussian, iff

$$
\forall T \in (0, \ \infty) \ \sup_{t \in [ 0, T]} ||X(t)||_{\psi(2)} < \infty.
$$

 \vspace{4mm}

 {\bf Theorem 4.1.} \ Let a deterministic variable  $ \ T \ $ be sufficiently great: $ \ T > e. \ $
   Suppose in addition to the conditions of theorem 2.1 that the "initial value" for $ \ \theta(\cdot), \ $  say, $ \ \theta(0) \ $ is subgaussian;
 for instance, may be a deterministic constant.  The described in definition 4.1 (\ref{inf operator Qb}) random process $ \ \theta = \theta(t), \ t \ge 0 \ $ is local subgaussian
 and has alike the Levi's Brownian motion  module of continuity:

\begin{equation} \label{module cont}
\exists \nu = \nu(\omega) \in G\psi(2), \ |\theta(t) - \theta(s)| \le \nu \cdot \sqrt{ |t - s| \ \ln(2 + 1/|t-s|)  }, \ t,s \in (0,1).
\end{equation}

 \ Further, there exists a subgaussian positive r.v. $ \ \tau, \ ||\tau||_{\psi(2)} = 1  \ $ such that

\begin{equation} \label{Log Law}
\max_{t \in [0,T]} |\theta(t)| \le C_1(\beta) \ \tau \ \sqrt{ T \ \ln T}.
\end{equation}

\vspace{4mm}

  \ {\bf Proof.} Suppose for beginning that $ \ ||\theta(s)||_{\psi(2)} < \infty. \ $ We apply the obtained estimate
(\ref{estim for p}) \ and obtain after simple calculations

\begin{equation} \label{distanse}
||\theta(t) - \theta(s)||_{\psi(2)} \le C_2(\beta) \ \sqrt{|t-s|}.
\end{equation}

\vspace{3mm}

 \ It remains to use the results of the sections 3.12 and 3.14 of the monograph \cite{Ostrovsky1}, devoting namely
the subgaussian random fields. \par

 \ The results (\ref{module cont}) and (\ref{Log Law}) are essentially non-improvable, for example, for the classical Bessel's processes
(\ref{Bessel n}). \par

\vspace{5mm}

\vspace{0.5cm} \emph{Acknowledgement.} {\footnotesize The first
author has been partially supported by the Gruppo Nazionale per
l'Analisi Matematica, la Probabilit\`a e le loro Applicazioni
(GNAMPA) of the Istituto Nazionale di Alta Matematica (INdAM) and by
Universit\`a degli Studi di Napoli Parthenope through the project
\lq\lq sostegno alla Ricerca individuale\rq\rq (triennio 2015-2017)}.\par

\vspace{4mm}

 \ The second author is very grateful to Molchanov S.A. for the statement of the considered problem and
valuable guidance. \par

\vspace{5mm}

\end{document}